\documentclass[12pt]{article}
\usepackage{amssymb}
\usepackage{amsmath}
\usepackage{amsthm}
\usepackage{graphicx}

\begin{document}

\title{A Variant of Multi-task n-vehicle Exploration Problem: Maximizing Every Processor's Average Profit\footnote
{Supported by 973 Program of China under Grant No. 2006CB701306, Key
Laboratory of Management, Decision and Information Systems, CAS, and
Beijing Research Center of Urban System Engineering.}}

\author{Yang-yang Xu\thanks{Email:xuyangyang@amss.ac.cn}, Jin-chuan
  Cui\thanks{Email:cjc@amss.ac.cn}}
\date{}

\maketitle

\begin{abstract}
We discuss a variant of multi-task n-vehicle exploration problem.
Instead of requiring an optimal permutation of vehicles in every
group, the new problem asks all vehicles in a group to arrive at a
same destination. It can also be viewed as to maximize every
processor's average profit, given $n$ tasks, and each task's
consume-time and profit. Meanwhile, we propose a new kind of
partition problem in fractional form, and analyze its computational
complexity. Moreover, by regarding fractional partition as a special
case, we prove that the maximizing average profit problem is NP-hard
when the number of processors is fixed and it is strongly NP-hard in
general. At last, a pseudo-polynomial time algorithm for the
maximizing average profit problem and the fractional partition
problem is presented, thanks to the idea of the pseudo-polynomial
time algorithm for the classical partition problem.
\end{abstract}



\section{Introduction}\label{Introduction}
Let $A=\{a_i\in \mathbf{Z}^+: i=1,\cdots,n\}$, and $B=\{b_i\in
\mathbf{Z}^+: i=1,\cdots,n\}$. The n-vehicle exploration problem was
introduced in [15] and analyzed in [8, 9, 13] to solve a scheduling
problem
$$\underset
{\pi}{\text{max}}\ \frac{a_{\pi(1)}}{b_{\pi(1)}}+\frac{a_{\pi(2)}}
{b_{\pi(1)}+b_{\pi(2)}}+\cdots+\frac{a_{\pi(n)}}{b_{\pi(1)}+\cdots+b_{\pi(n)}},$$
where $\left(\pi(1),\pi(2),\cdots,\pi(n)\right)$ is a permutation of
$(1,2,\cdots,n)$. It can be described as$^{[15]}$: given $n$
vehicles, the $i$th one can carry at most $a_i$ liters of oil, and
consumes $b_i$ liters of oil per kilometer, $i=1,\cdots,n$. These
vehicles start to travel towards one direction from the same
position at the same time. On the path, they can not get oil from
outside, but at any position any one can stop and transfer its left
oil to other vehicles. How can we arrange these vehicles to make one
of them travel farthest and ensure that all of them return to the
original position?

The authors of [14] extended the n-vehicle exploration problem into
a multi-task n-vehicle exploration problem (MTNVEP). They first
divide $n$ vehicles into $m$ groups and then solve $m$ new n-vehicle
exploration problem to make every group travel far enough. That is
to solve
\begin{eqnarray}
\text{max}\underset{1\le j\le m}{\text{min}}\underset
{\pi_j}{\text{max}}\
\frac{a_{\pi_j(1)}}{b_{\pi_j(1)}}+\frac{a_{\pi_j(2)}}
{b_{\pi_j(1)}+b_{\pi_j(2)}}+\cdots+\frac{a_{\pi_j(n_j)}}{b_{\pi_j(1)}+\cdots+b_{\pi_j(n_j)}},
\end{eqnarray} where
$(1_1,\cdots,1_{n_1}),\cdots,(m_1,\cdots,m_{n_m})$ is a partition of
$(1,2,\cdots,n)$, and $\left(\pi_j(1),\cdots,\pi_j(n_j)\right)$ is a
permutation of $(j_1,\cdots,j_{n_j}), j=1,\cdots,m$. Note that the
first ``max'' objective of (1) is about $(1,2,\cdots,n)$'s
partition. They concluded that MTNVEP is NP-hard when $m$ is fixed,
and strongly NP-hard for general $m$. We notice that in their proof
MTNVEP's computational complexity is not related to vehicles'
permutation. Namely, even if the n-vehicle exploration problem can
be solved in a polynomial time, MTNVEP is still NP-hard.

In this paper, we drop the permutation requirement of MTNVEP, but
ask all vehicles in every group to arrive at a same destination.
That is to solve
\begin{eqnarray}
\text{max}\underset{1\le j\le
m}{\text{min}}\frac{a_{j_1}+\cdots+a_{j_{n_j}}}{b_{j_1}+\cdots+b_{j_{n_j}}},
\end{eqnarray}
where $(1_1,\cdots,1_{n_1}),\cdots,(m_1,\cdots,m_{n_m})$ is a
partition of $(1,2,\cdots,n)$. The ``max'' objective of (2) is still
about $(1,\cdots,n)$'s partition. This is a variant of multi-task
n-vehicle exploration problem, which can be viewed as to maximize
every processor's average profit (MAP). It can be described as:
given $n$ tasks, and $m$ identical processors, it will take $b_i$
units of time to finish $i$th task and make $a_i$ units of profit.
How can we distribute the $n$ tasks to $m$ processors so that every
processor can make its average profit large enough?

Meanwhile, we define a new kind of partition problem in fractional
form (FP), if we change our objective to determine whether there
exists a partition of $(1,2,\cdots,n)$ into
$(1_1,\cdots,1_{n_1}),\cdots,(m_1,\cdots,m_{n_m})$ such that
$$\frac{a_{j_1}+\cdots+a_{j_{n_j}}}{b_{j_1}+\cdots+b_{j_{n_j}}}
=\frac{a_1+a_2+\cdots+a_n}{b_1+b_2+\cdots+b_n},\text{ for any
}j=1,\cdots,m,$$ instead of solving (2). FP is very similar to the
classical partition problem proposed by Karp in his twenty-one
NP-complete problems$^{[6]}$, and 3-partition problem proposed by
Garey and Johnson$^{[2, 3]}$. The classical partition problem
locates at the center of computationally intractable problems. Based
on the partition problem, many other hard problems are proved to be
NP-complete or NP-hard, such as bin packing$^{[4]}$, multiprocessor
scheduling problem$^{[5]}$, and 0-1 integer programming$^{[11]}$.
3-Partition is often employed to devise a strong NP-completeness
proof. We will prove FP's NP-completeness by reducing partition
problem to FP when $m$ is fixed, and prove FP's strong
NP-completeness by reducing 3-partition to FP for general $m$.
Meanwhile, we will prove that MAP is NP-hard$^{[3, 12]}$ for fixed
$m\ge2$, and strongly NP-hard for general $m\ge2$ regarding FP as a
special case of MAP. Moreover, thanks to the idea of designing
pseudo-polynomial time algorithm for classical partition
problem$^{[3, 7]}$, we design a pseudo-polynomial time algorithm for
MAP and FP.

The outline of this paper is as follows: In Section 2, we give a
general description of MAP and FP. Their complexity is analyzed in
Section 3. We prove that FP is NP-complete and MAP is NP-hard for a
fixed $m\ge2$, and that FP is strongly NP-complete and MAP is
strongly NP-hard for general $m\ge2$. In Section 4, a
pseudo-polynomial time algorithm for MAP and FP is presented, thanks
to the idea of the pseudo-polynomial time algorithm for the
partition problem. Section 5 concludes the paper.
\section{The models}
\newtheorem{definition}{Definition}
\begin{definition}[Maximizing average profit]
Given $A=\{a_i\in \mathbf{Z}^+: i=1,\cdots,n\}$, and $B=\{b_i\in
\mathbf{Z}^+: i=1,\cdots,n\}$, partition the index set
$I=\{1,2,\cdots,n\}$ into $m$ disjoint subsets $I_1,I_2,\cdots,I_m$,
denote
$$r_j=\frac{\underset{i\in I_j}{\sum} a_i}{\underset{i\in I_j}{\sum} b_i}\ ,\quad j=1,\cdots,m,$$ and let
$$f(I_1,I_2,\cdots,I_m)=\underset{1\le j\le m}{min}r_j.$$
Determine the maximum value of $f$.
\end{definition}

\begin{definition}[Fractional partition]
Given $A=\{a_i\in \mathbf{Z}^+: i=1,\cdots,n\}$, $B=\{b_i\in
\mathbf{Z}^+: i=1,\cdots,n\}$, $I=\{1,2,\cdots,n\}$, and
$S=\underset{i\in I}{\sum}a_i$, $T=\underset{i\in I}{\sum}b_i$,
partition $I$ into $m$ disjoint subsets $I_1,I_2,\cdots,I_m$, and
denote
$$r_j=\frac{\underset{i\in I_j}{\sum} a_i}{\underset{i\in I_j}{\sum} b_i}\ ,\quad j=1,\cdots,m.$$
Determine whether there exists a partition of $I$ such that
$$r_j=\frac{S}{T},\ j=1,\cdots,m.$$
\end{definition}

\newtheorem{remark}{Remark}

\section{Complexity analysis}\label{complexity}
In the first subsection, we will prove that: for fixed $m\ge2$, FP
is NP-complete and MAP is NP-hard; in the second subsection, we will
prove that: for general $m\ge 2$, FP is strongly NP-complete and MAP
is strongly NP-hard.
\subsection{NP-completeness}
\newtheorem{lemma}{Lemma}

\newtheorem{theorem}{Theorem}
\begin{theorem}
FP is NP-complete for fixed $m\ge2$.
\end{theorem}

\begin{proof}
At first, we prove this theorem for $m=2$ in three steps.

{\it Step 1.} It is easy to see that FP belongs to NP, because given
a partition of $I$, we can immediately calculate $r_j$ and check
whether $r_j={S}/{T}$, for $j=1,\cdots,m$.

{\it Step 2.} We reduce one instance of partition problem to FP in a
polynomial time.

The instance of partition problem is:

$Q_1$: Given $C=\{c_i\in \mathbf{Z}^+,\ i=1,\cdots,n\}$, and
$2K=\underset{i}{\sum} c_i$, determine whether there exists a subset
$C_1\subseteq C$, such that
    $$\sum_{c\in C_1}^{} c=\sum_{c\in C-C_1}^{} c=K.$$

We construct FP's instance as follows:

$Q_2$: Let
\begin{eqnarray*}
&&A=\left\{Mc_1,\cdots,Mc_n,
\underset{N+n-1}{\underbrace{M\delta,\cdots,M\delta}},\underset{N-2n+1}{\underbrace{3M\delta/2,\cdots,3M\delta/2}},MK,MK\right\}\\
&&B=\left\{\underset{2N}{\underbrace{M,\cdots,M}},MN+M\epsilon,MN+M\epsilon
\right\}
\end{eqnarray*}
where
\begin{eqnarray*}
&&N=4K+1,\ \epsilon=1/\lceil1/\epsilon'\rceil,\
M=(5N-4n+1)/\epsilon,\
\delta=\dfrac{2\epsilon}{5N-4n+1}\\
&&\epsilon'=\dfrac{-(4N^3+2KN-N^2-1)+\sqrt{(4N^3+2KN-N^2-1)^2+16N^3}}{8N^2};
\end{eqnarray*}
Determine whether there exists a partition of the index set
$I$($:=\{1,2,\cdots,2N+2\}$) into two disjoint subsets $I_1,I_2$,
such that
$$r_1=r_2=\frac{2K+\epsilon/2}{2N+\epsilon}.$$

It is easy to see that the reduction can be finished in a polynomial
time. Though $N,M$, and $K$ are not bounded by any polynomial of
$n$, all numbers can be obtained in polynomial time. Moreover, we
need construct $M\delta$ and $3M\delta/2$ only once, and then denote
the number of times they appear in $A$ by $N+n-1$ and $N-2n+1$,
respectively.

Note that all numbers in $A$ and $B$ have a common divisor $M$. This
is to satisfy the integer requirement. Since the numerator and
denominator of every fraction will have a common divisor $M$, it
will not change the result to eliminate $M$ in every number at the
beginning. Thus we will consider the following sets
\begin{eqnarray*}
&&A=\left\{c_1,\cdots,c_n,
\underset{N+n-1}{\underbrace{\delta,\cdots,\delta}},\underset{N-2n+1}{\underbrace{3\delta/2,\cdots,3\delta/2}},K,K\right\}\\
&&B=\left\{\underset{2N}{\underbrace{1,\cdots,1}},N+\epsilon,N+\epsilon
\right\}
\end{eqnarray*}

{\it Step 3.} We prove that $Q_1$ is true if and only if $Q_2$ is
true, where $Q_i$ ($i=1,2$) is true if there exists a partition as
described.

If $Q_1$ is true, i.e., there exists a subset $C_1\subseteq C$, such
that $$\underset{c\in C_1}{\sum}c=\underset{c\in C-C_1}{\sum}c=K,$$
denote $n_1=|C_1|, n_2=n-n_1$, and set
\begin{eqnarray*}
&&I_1=\{i:c_i\in C_1\}\cup\{n+1,\cdots,\frac{N-1}{2}+3n_2,
N+2n,\cdots,\frac{3N-1}{2}+2n_1,2N+1\}\\
&&I_2=I-I_1.
\end{eqnarray*}
It is not difficult to verify that
$r_1=r_2=(2K+\epsilon/2)/(2N+\epsilon)$. Namely, $Q_2$ is true.

Conversely, if $Q_2$ is true, i.e., there exists $I$'s partition
$I_1,I_2$, such that
\begin{eqnarray}
r_1=r_2=\frac{2K+\epsilon/2}{2N+\epsilon},
\end{eqnarray}
assume $\underset{i\in I_j}{\sum}c_i=p_j$, the numerator of fraction
$r_j$ is $p_j+y_j\delta+x_jK$, and the denominator of $r_j$ is
$q_j+x_j(N+\epsilon)$, where $x_j, y_j, q_j(j=1,2)$ are nonnegative
integers. Then we get
\begin{eqnarray}
r_j=\frac{p_j+y_j\delta+x_jK}{q_j+x_j(N+\epsilon)}=\frac{2K+\epsilon/2}{2N+\epsilon},\
j=1,2.
\end{eqnarray}
Note that $q_j\neq0$; otherwise $p_j=y_j=0$, so
$r_j=K/(N+\epsilon)$. A contradiction to (3). Next we will prove
$$p_j=K,\ j=1,2.$$

First, it must be $p_j/q_j\ge K/N,\ j=1,2$; otherwise, without loss
of generality, we assume $p_1/q_1<K/N$, then we have (because
$q_1\le 2N$)
\begin{eqnarray}
\frac{p_1}{q_1}\le\frac{K}{N}-\frac{1}{2N^2}.
\end{eqnarray}
Transforming (4), we obtain
\begin{eqnarray}
\frac{p_1}{q_1}=\frac{2K+\epsilon/2}{2N+\epsilon}+\frac{Nx_1\epsilon/2
+x_1\epsilon(K+\epsilon/2)}{q_1(2N+\epsilon)}-\frac{\epsilon
y_1\delta+2Ny_1\delta}{q_1(2N+\epsilon)}.
\end{eqnarray}
Using $N=4K+1\ge 2n$ and $\delta=2\epsilon/(5N-4n+1)$, it follows
from (5) and (6) that
\begin{eqnarray}
\frac{K}{N}-\frac{1}{2N^2}\ge\frac{p_1}{q_1}\ge\frac{2K+\epsilon/2}{2N+\epsilon}-\epsilon.
\end{eqnarray}
We can choose a sufficiently small $\epsilon>0$ such that
$$\frac{K}{N}-\frac{1}{2N^2}<\frac{2K+\epsilon/2}{2N+\epsilon}-\epsilon.$$
That is
$$2N^2\epsilon^2+(4N^3+2KN-N^2-1)\epsilon-2N<0.$$
If we choose $\epsilon'$ as half of the positive root of equation
$$2N^2\epsilon^2+(4N^3+2KN-N^2-1)\epsilon-2N=0,$$
i.e.
$$\epsilon'=\dfrac{-(4N^3+2KN-N^2-1)+\sqrt{(4N^3+2KN-N^2-1)^2+16N^3}}{8N^2},$$
and let $\epsilon=1/\lceil1/\epsilon'\rceil$,
then$$\frac{K}{N}-\frac{1}{2N^2}<\frac{2K+\epsilon/2}{2N+\epsilon}-\epsilon.$$
A contradiction to (7). Thus $p_1/q_1\ge K/N$, and $p_2/q_2\ge K/N$.

Secondly, because $(p_1+p_2)/(q_1+q_2)=K/N$, we get
\begin{eqnarray}
\frac{p_j}{q_j}=\frac{K}{N},\ j=1,2.
\end{eqnarray}
Write (8) as $2p_jN=2q_jK$. Note that $N$ and $2K$ are relatively
prime$^{[1]}$, so $q_j$ is divisible by $N$. That is, $q_j=N$ or
$2N$, but $q_1+q_2=2N$, thus $q_j=N,\ j=1,2$. Together with (8), we
have
$$p_j=K,j=1,2.$$
If we set $$C_1=\{c_i\in C|i\in I_1\},$$ then
$$\underset{c\in C_1}{\sum}c=\underset{c\in C-C_1}{\sum}c=p_1=K.$$
Namely, $Q_1$ is true.

Until now, we have proved that FP is NP-complete when $m=2$; for
$m>2$, we can construct FP's instance as follows.

$Q_2^{'}$: Let
\begin{eqnarray*}
&&A=\left\{Mc_1,\cdots,Mc_n,
\underset{N+n-1}{\underbrace{M\delta,\cdots,M\delta}},\underset{N-2n+1}{\underbrace{3M\delta/2}},MK,MK,
\underset{m-2}{\underbrace{2MK+M\epsilon/2}}\right\}\\
&&B=\left\{\underset{2N}{\underbrace{M,\cdots,M}},MN+M\epsilon,MN+M\epsilon,\underset{m-2}{\underbrace{2MN+M\epsilon}}
\right\}
\end{eqnarray*}
where
\begin{eqnarray*}
&&N=4K+1,\ \epsilon=1/\lceil1/\epsilon'\rceil,\
M=(5N-4n+1)/\epsilon,\
\delta=\dfrac{2\epsilon}{5N-4n+1}\\
&&\epsilon'=\dfrac{-(4N^3+2KN-N^2-1)+\sqrt{(4N^3+2KN-N^2-1)^2+16N^3}}{8N^2};
\end{eqnarray*}
Determine whether there exists a partition of the index set
$I$($:=\{1,2,\cdots,2N+m\}$) into $m$ disjoint subsets
$I_1,I_2,\cdots,I_m$, such that
$$r_j=\frac{2K+\epsilon/2}{2N+\epsilon},\ j=1,\cdots,m.$$

This construction can also be finished in polynomial time. The proof
of equivalence between $Q_1$ and $Q_2^{'}$ is essentially the same
as that when $m=2$.

Therefore, FP is NP-complete for fixed $m\ge2$.
\end{proof}

\begin{theorem}
MAP is NP-hard for fixed $m\ge2$.
\end{theorem}

\begin{proof}
Given $A=\{a_i\in \mathbf{Z}^+: i=1,\cdots,n\}$, $B=\{b_i\in
\mathbf{Z}^+: i=1,\cdots,n\}$, $I=\{1,2,\cdots,n\}$, and
$S=\underset{i\in I}{\sum}a_i$, $T=\underset{i\in I}{\sum}b_i$,
partition $I$ into $m$ disjoint subsets $I_1,I_2,\cdots,I_m$, and
denote
$$r_j=\frac{\underset{i\in I_j}{\sum} a_i}{\underset{i\in I_j}{\sum} b_i}\ ,\quad j=1,\cdots,m.$$
Then the proof of this theorem is straightforward from ``$r_j\ge
S/T$ for any $j$ if and only if $r_j=S/T$ for any $j$.''
\end{proof}

\newtheorem{corollary}{Corollary}

\subsection{Strong NP-completeness}

\begin{theorem}
FP is strongly NP-complete for general $m\ge2$.
\end{theorem}

\begin{proof}
In last subsection, we have shown that FP belongs to NP. Next, we
only need reduce 3-Partition to FP for general $m\ge2$ in
pseudo-polynomial time.

An instance of 3-Partition is:

$Q_3$: Given $D=\{d_i\in \mathbf{Z}^+:i=1,\cdots,3m\},\
mK=\underset{i}{\sum}d_i$, and $K/4<d_i<K/2$ for every $i$,
determine whether there exists a partition of $D$ into $m$ disjoint
subsets $D_1,\cdots,D_m$, such that $$\underset{d\in D_j}{\sum}d=K,\
j=1,\cdots,m.$$

We construct FP's instance as follows:

$Q_4$: Let
\begin{eqnarray*}
&&A=\left\{Md_1,\cdots,Md_{3m},\underset{mN-3m}{\underbrace{M\delta}},\underset{m}{\underbrace{MK}}\right\}\\
&&B=\left\{\underset{mN}{\underbrace{M,\cdots,M}},\underset{m}{\underbrace{MN+M\epsilon}}\right\}
\end{eqnarray*}
where
\begin{eqnarray*}
&&N=2mK+1,\ \epsilon=1/\lceil1/\epsilon'\rceil,\
M=(mN-3m)/\epsilon,\
\delta=\dfrac{\epsilon}{mN-3m}\\
&&\epsilon'=\dfrac{-(2mN^3+mKN-1)+\sqrt{(2mN^3+mKN-1)^2+8mN^3}}{4mN^2};
\end{eqnarray*}
Determine whether there exists a partition of the index set
$I$($:=\{1,2,\cdots,m(N+1)\}$) into $m$ disjoint subsets
$I_1,I_2,\cdots,I_m$, such that
$$r_j=\frac{2K+\epsilon/m}{2N+\epsilon},\ j=1,\cdots,m.$$

We bring a common divisor $M$ in each element of $A$ and $B$, for
the same reason as we proved theorem 1. It will not change the
result to eliminate $M$ at the beginning. Thus we will consider the
following sets
$$A=\left\{d_1,\cdots,d_{3m},\underset{mN-3m}{\underbrace{\delta}},\underset{m}{\underbrace{K}}\right\},\quad
B=\left\{\underset{mN}{\underbrace{1,\cdots,1}},\underset{m}{\underbrace{N+\epsilon}}\right\}.$$

Now we focus on proving that $Q_3$ is true if and only if $Q_4$ is
true.

If $Q_3$ is true, i.e., there exists a partition of $D$ into $m$
disjoint subsets $D_1,\cdots,D_m$, such that $$\underset{d\in
D_j}{\sum}d=K,\ j=1,\cdots,m,$$ there are exactly three elements in
every subset, since $K/4<d_i<K/2,\ i=1,\cdots,3m$. Letting
$$I_j=\left\{i:d_i\in D_j\right\}\cup\{3m+(j-1)(N-3)+1,\cdots,3m+j(N-3),mN+j\},$$
we have
$$r_j=\frac{2K+(N-3)\delta}{2N+\epsilon}=\frac{2K+\epsilon/m}{2N+\epsilon},\ j=1,\cdots,m.$$
Namely, $Q_4$ is true.

Conversely, if $Q_4$ is true, i.e., there exists a partition of the
index set $I$ into $m$ disjoint subsets $I_1,\cdots,I_m$, such that
\begin{eqnarray}
r_j=\frac{2K+\epsilon/m}{2N+\epsilon},\ j=1,\cdots,m,
\end{eqnarray}
assume $\underset{i\in I_j}{\sum}d_i=p_j$, the numerator of $r_j$ is
$p_j+x_jK+y_j\delta$, and the denominator of $r_j$ is
$q_j+x_j(N+\epsilon)$, , where $x_j, y_j, q_j(j=1,2)$ are
nonnegative integers. Then we have
$$r_j=\frac{p_j+x_jK+y_j\delta}{q_j+x_j(N+\epsilon)},=\frac{2K+\epsilon/m}{2N+\epsilon},\ j=1,\cdots,m.$$
Note that $q_j\ge 1$ for every $j=1,\cdots,m$; otherwise
$p_j=y_j=0$, then $r_j={K}/{(N+\epsilon)}$. A contradiction to (9).
Next, we will prove $$p_j=K,\ j=1,\cdots,m$$

Firstly, it must be $p_j/q_j\ge K/N, j=1,\cdots,m$. Otherwise,
without loss of generality, we assume $p_1/q_1< K/N$, then we have
(because $q_1\le mN$)
$$\frac{p_1}{q_1}\le \frac{K}{N}-\frac{1}{mN^2}.$$
Because
$$r_1=\frac{p_1+x_1K+y_1\delta}{q_1+x_1(N+\epsilon)}\le\frac{p_1+x_1K+\epsilon}{q_1+x_1(N+\epsilon)}$$
and $$r_1=\frac{2K+\epsilon/m}{2N+\epsilon}>
\frac{2K}{2N+\epsilon},$$ we have
\begin{eqnarray}\frac{p_1+x_1K+\epsilon}{q_1+x_1(N+\epsilon)}>\frac{2K}{2N+\epsilon}.
\end{eqnarray}
Transforming (10), we get
\begin{eqnarray}
\frac{p_1}{q_1}&>&\frac{2K}{2N+\epsilon}+\frac{\epsilon
x_1K}{q_1(2N+\epsilon)}-\frac{\epsilon}{q_1}\nonumber\\
&\ge&\frac{2K}{2N+\epsilon}+\frac{\epsilon
x_1K}{q_1(2N+\epsilon)}-\epsilon\nonumber\\
&\ge&\frac{2K}{2N+\epsilon}-\epsilon.
\end{eqnarray}
We can choose a sufficiently small $\epsilon>0$, such that
$$\frac{K}{N}-\frac{1}{mN^2}<\frac{2K}{2N+\epsilon}-\epsilon;$$
that is $$mN^2\epsilon^2+(2mN^3+mKN-1)\epsilon-2N<0.$$ If we choose
$\epsilon'$ as half of the positive root of equation
$$mN^2\epsilon^2+(2mN^3+mKN-1)\epsilon-2N=0,$$
i.e.
$$\epsilon'=\dfrac{-(2mN^3+mKN-1)+\sqrt{(2mN^3+mKN-1)^2+8mN^3}}{4mN^2},$$
and let $\epsilon=1/\lceil1/\epsilon'\rceil$, then
$$\frac{K}{N}-\frac{1}{mN^2}<\frac{2K}{2N+\epsilon}-\epsilon.$$
But$$\frac{p_1}{q_1}\le\frac{K}{N}-\frac{1}{mN^2},$$ a contradiction
to (11). Thus, for the chosen $\epsilon$, it must be
$$\frac{p_j}{q_j}\ge\frac{K}{N},\ j=1,\cdots,m.$$

Secondly, because $$\frac{p_j}{q_j}\ge\frac{K}{N},\ j=1,\cdots,m,$$
and
$$\frac{\sum p_j}{\sum q_j}=\frac{mK}{mN}=\frac{K}{N},$$
we get
\begin{eqnarray}
\frac{p_j}{q_j}=\frac{K}{N},\ j=1,\cdots,m.
\end{eqnarray}
Write (12) as $mp_jN=mq_jK,\ j=1,\cdots,m$. Note that $N$ and $mK$
are relatively prime, thus $q_j$ is divisible by $N$. Because
$\sum_{j=1}^m q_j=mN$, we get $q_j=N,\ j=1,\cdots,m$. Together with
(12), we have
$$p_j=K,\ j=1,\cdots,m.$$ If we set
$$D_j=\{d_i\in D:i\in I_j\},\ j=1,\cdots,m,$$ then
$$\underset{d\in D_j}{\sum}d=p_j=K,\ j=1,\cdots,m.$$ Namely, $Q_3$ is
true.

Therefore, FP is strongly NP-complete for general $m\ge2$.
\end{proof}

\begin{theorem}
MAP is strongly NP-hard for general $m\ge2$.
\end{theorem}
\begin{proof}
The proof of this theorem is essentially the same as that of theorem
2.
\end{proof}

\section{Pseudo-polynomial algorithm}
In last section, we analyzed the computational complexity of FP and
MAP: FP is NP-complete and MAP is NP-hard for fixed $m\ge2$; FP is
strongly NP-complete and MAP is strongly NP-hard for general
$m\ge2$. According to [3], if problem $\prod$ is strongly
NP-complete, then it is impossible to solve $\prod$ in a
pseudo-polynomial time unless P=NP. As a result, we do not hope to
find FP and MAP's pseudo-polynomial time algorithm for general
$m\ge2$. However, we can prove that FP and MAP can be solved in a
pseudo-polynomial time when $m$ is fixed. The pseudo-polynomial time
algorithm we will present for FP and MAP is similar to that for the
classical partition problem.

\begin{theorem}
FP can be solved in a pseudo-polynomial time for fixed $m\ge2$.
\end{theorem}

\begin{proof}
To prove the theorem, we design a pseudo-polynomial time algorithm
of FP when $m$ is fixed.

Define a boolean function
$t(i,{p_1}/{q_1},\cdots,{p_{m-1}}/{q_{m-1}})$, meaning whether there
exist $m-1$ disjoint subsets $I_1,\cdots,I_{m-1}$ of
$\{1,2,\cdots,i\}$, such that the $m-1$ fractions $(\sum_{k\in
I_1}a_k)/(\sum_{k\in I_1}b_k)$, $\cdots$, $(\sum_{k\in
I_{m-1}}a_k)/(\sum_{k\in I_{m-1}}b_k)$ are exactly
$p_1/q_1,\cdots,{p_{m-1}}/{q_{m-1}}$, respectively. Note that we ask
the fractions to be same not only in value, but also in form. For
example, ${1}/{1}$ is not considered the same as ${2}/{2}$.
Moreover, we prescribe that the fraction $(\sum_{k\in
I_j}a_k)/(\sum_{k\in I_j}b_k)$ is ${0}/{0}$ if $I_j=\emptyset$, and
the value of fraction ${p}/{0}$ is {\it zero}. If there exist
$I_1,\cdots,I_{m-1}$ as described, then
$$t(i,\frac{p_1}{q_1},\cdots,\frac{p_{m-1}}{q_{m-1}})=1;$$ else
$$t(i,\frac{p_1}{q_1},\cdots,\frac{p_{m-1}}{q_{m-1}})=0.$$

It is easy to verify that
$$t(i+1,\frac{p_1}{q_1},\cdots,\frac{p_{m-1}}{q_{m-1}})=1,$$ if and
only if $$t(i,\frac{p_1}{q_1},\cdots,\frac{p_{m-1}}{q_{m-1}})=1$$ or
$$t(i,\frac{p_1}{q_1},\cdots,\frac{p_k-a_{i+1}}{q_k-b_{i+1}},\cdots,\frac{p_{m-1}}{q_{m-1}})=1$$
for some $k$ with $1\le k\le m-1$. Then we can list the following
matrices in order

\begin{eqnarray*}
\left[\begin{array}{ccc}
t(1,\frac{0}{0},\frac{0}{0},\cdots,\frac{0}{0}),t(1,\frac{0}{1},\frac{0}{0},\cdots,\frac{0}{0}),
\cdots,t(1,\frac{0}{T},\frac{0}{0},\cdots,\frac{0}{0})\\[0.2cm]
t(1,\frac{1}{0},\frac{0}{0},\cdots,\frac{0}{0}),t(1,\frac{1}{1},\frac{0}{0},\cdots,\frac{0}{0}),
\cdots,t(1,\frac{1}{T},\frac{0}{0},\cdots,\frac{0}{0})\\[0.2cm]
\cdots\cdots\cdots\cdots\cdots\cdots\cdots\cdots\cdots\cdots\cdots\cdots\cdots\cdots\cdots\\[0.2cm]
t(1,\frac{S}{0},\frac{0}{0},\cdots,\frac{0}{0}),t(1,\frac{S}{1},\frac{0}{0},\cdots,\frac{0}{0}),
\cdots,t(1,\frac{S}{T},\frac{0}{0},\cdots,\frac{0}{0})
\end{array}
\right]\cdots
\end{eqnarray*}

\begin{eqnarray*}
\left[\begin{array}{ccc}
t(1,\frac{0}{0},\frac{p_2}{q_2},\cdots,\frac{p_{m-1}}{q_{m-1}}),
t(1,\frac{0}{1},\frac{p_2}{q_2},\cdots,\frac{p_{m-1}}{q_{m-1}}),
\cdots,t(1,\frac{0}{T},\frac{p_2}{q_2},\cdots,\frac{p_{m-1}}{q_{m-1}})\\[0.2cm]
t(1,\frac{1}{0},\frac{p_2}{q_2},\cdots,\frac{p_{m-1}}{q_{m-1}}),
t(1,\frac{1}{1},\frac{p_2}{q_2},\cdots,\frac{p_{m-1}}{q_{m-1}})
\cdots,t(1,\frac{1}{T},\frac{p_2}{q_2},\cdots,\frac{p_{m-1}}{q_{m-1}})\\[0.2cm]
\cdots\cdots\cdots\cdots\cdots\cdots\cdots\cdots\cdots\cdots
\cdots\cdots\cdots\cdots\cdots\cdots\cdots\cdots\cdots\cdots\\[0.2cm]
t(1,\frac{S}{0},\frac{p_2}{q_2},\cdots,\frac{p_{m-1}}{q_{m-1}}),
t(1,\frac{S}{1},\frac{p_2}{q_2},\cdots,\frac{p_{m-1}}{q_{m-1}}),
\cdots,t(1,\frac{S}{T},\frac{p_2}{q_2},\cdots,\frac{p_{m-1}}{q_{m-1}})
\end{array}
\right]\cdots
\end{eqnarray*}

\begin{eqnarray*}
\left[
\begin{array}{ccc}
t(1,\frac{S}{T},\cdots,\frac{S}{T},\frac{0}{0}),
t(1,\frac{S}{T},\cdots,\frac{S}{T},\frac{0}{1}),
\cdots,t(1,\frac{S}{T},\cdots,\frac{S}{T},\frac{0}{T})\\[0.2cm]
t(1,\frac{S}{T},\cdots,\frac{S}{T},\frac{1}{0}),
t(1,\frac{S}{T},\cdots,\frac{S}{T},\frac{1}{1}),
\cdots,t(1,\frac{S}{T},\cdots,\frac{S}{T},\frac{1}{T})\\[0.2cm]
\cdots\cdots\cdots\cdots\cdots\cdots
\cdots\cdots\cdots\cdots\cdots\cdots\cdots\cdots\cdots\cdots\\[0.2cm]
t(1,\frac{S}{T},\cdots,\frac{S}{T},\frac{S}{0}),
t(1,\frac{S}{T},\cdots,\frac{S}{T},\frac{S}{1}),
\cdots,t(1,\frac{S}{T},\cdots,\frac{S}{T},\frac{S}{T})
\end{array}
\right]\cdots
\end{eqnarray*}

\begin{eqnarray*}
\left[
\begin{array}{ccc}
t(n-1,\frac{S}{T},\cdots,\frac{S}{T},\frac{0}{0}),
t(n-1,\frac{S}{T},\cdots,\frac{S}{T},\frac{0}{1}),
\cdots,t(n-1,\frac{S}{T},\cdots,\frac{S}{T},\frac{0}{T})\\[0.2cm]
t(n-1,\frac{S}{T},\cdots,\frac{S}{T},\frac{1}{0}),
t(n-1,\frac{S}{T},\cdots,\frac{S}{T},\frac{1}{1}),
\cdots,t(n-1,\frac{S}{T},\cdots,\frac{S}{T},\frac{1}{T})\\[0.2cm]
\cdots\cdots\cdots\cdots\cdots\cdots\cdots\cdots
\cdots\cdots\cdots\cdots\cdots\cdots\cdots\cdots\cdots\cdots\\[0.2cm]
t(n-1,\frac{S}{T},\cdots,\frac{S}{T},\frac{S}{0}),
t(n-1,\frac{S}{T},\cdots,\frac{S}{T},\frac{S}{1}),
\cdots,t(n-1,\frac{S}{T},\cdots,\frac{S}{T},\frac{S}{T})
\end{array}
\right]
\end{eqnarray*}
While listing these matrices, if some
$$t(i,\frac{p_1}{q_1},\cdots,\frac{p_{m-1}}{q_{m-1}})=1,$$ and
\begin{eqnarray*}
&&\frac{p_j}{q_j}=\frac{S}{T},\ j=1,\cdots,m-1\\
&&\frac{S-\overset{m-1}{\underset{j=1}{\sum}}p_j}{T-\overset{m-1}{\underset{j=1}{\sum}}q_j}=\frac{S}{T},
\end{eqnarray*} the
algorithm stops and return ``FP is true''; else return ``FP is
false''.

For example, when $m=2$, we can list
\begin{eqnarray*}
&&\left[\begin{array}{ccc} t(1,\frac{0}{0}),t(1,\frac{0}{1}),
\cdots,t(1,\frac{0}{T})\\[0.2cm]
t(1,\frac{1}{0}),t(1,\frac{1}{1}),
\cdots,t(1,\frac{1}{T})\\[0.2cm]
\cdots\cdots\cdots\cdots\cdots\cdots\\[0.2cm]
t(1,\frac{S}{0}),t(1,\frac{S}{1}), \cdots,t(1,\frac{S}{T})
\end{array}
\right]\cdots \left[\begin{array}{ccc}
t(i,\frac{0}{0}),t(i,\frac{0}{1}),
\cdots,t(i,\frac{0}{T})\\[0.2cm]
t(i,\frac{1}{0}),t(i,\frac{1}{1}),
\cdots,t(i,\frac{1}{T})\\[0.2cm]
\cdots\cdots\cdots\cdots\cdots\cdots\\[0.2cm]
t(i,\frac{S}{0}),t(i,\frac{S}{1}), \cdots,t(i,\frac{S}{T})
\end{array}
\right]\\[0.5cm]
 &&\cdots\left[\begin{array}{ccc}
t(n-1,\frac{0}{0}),t(n-1,\frac{0}{1}),
\cdots,t(n-1,\frac{0}{T})\\[0.2cm]
t(n-1,\frac{1}{0}),t(n-1,\frac{1}{1}),
\cdots,t(n-1,\frac{1}{T})\\[0.2cm]
\cdots\cdots\cdots\cdots\cdots\cdots\\[0.2cm]
t(n-1,\frac{S}{0}),t(n-1,\frac{S}{1}), \cdots,t(n-1,\frac{S}{T})
\end{array}
\right]
\end{eqnarray*}

In all matrices, there are $(n-1)(ST)^{m-1}$ entries in total, and
to get the value of $t(i,{p_1}/{q_1},\cdots,{p_{m-1}}/{q_{m-1}})$
needs at most $m$ steps, so the computational complexity is
$O(mn(ST)^{m-1})$. Therefore, FP can be solved in a
pseudo-polynomial time.
\end{proof}

In the same way, we can solve MAP in a pseudo-polynomial time.

\begin{theorem}
MAP can be solved in a pseudo-polynomial time for fixed $m\ge2$.
\end{theorem}

\begin{proof}
Define
\begin{eqnarray*}
g(i,\frac{p_1}{q_1},\cdots,\frac{p_{m-1}}{q_{m-1}})=\left\{
\begin{array}{lll}
0&if\ \ t(\cdot)=0\\
min(\frac{p_1}{q_1},\cdots,\frac{p_{m-1}}{q_{m-1}},
\frac{S-\overset{m-1}{\underset{j=1}{\sum}}p_j}{T-\overset{m-1}{\underset{j=1}{\sum}}q_j}
)&if\ \ t(\cdot)=1
\end{array}
\right.
\end{eqnarray*}

While listing matrices in the proof of theorem 5, we record the
value of function {\it g} together with the boolean function {\it
t}. At the end, we search for all entries in these matrices to
obtain the maximum value of {\it g}, which is the optimal value of
MAP. To get the value of
$g(i,{p_1}/{q_1},\cdots,{p_{m-1}}/{q_{m-1}})$ needs at most $m$
steps, and to search for all entries needs $(n-1)(ST)^{m-1}$ steps,
so the computational complexity is $O(mn(ST)^{m-1})$. Therefore, MAP
can be solved in pseudo-polynomial time.
\end{proof}

\section{Conclusion and future work}
The multi-task n-vehicle exploration problem is very important in
the context of how to use limited resources to complete a number of
tasks. Due to the fact that MTNVEP is still NP-hard even if the
n-vehicle exploration problem can be solved in polynomial time, we
try to drop the permutation requirement but ask all vehicles in one
group to arrive at a same destination. The new problem is a variant
of multi-task n-vehicle exploration problem. It can also be regarded
as how to maximize every processor's average profit. Through
analyzing its special case: fractional partition, we prove that MAP
is NP-hard when the processor number is fixed, and it is strongly
NP-hard in general.

Fractional partition is a new kind of partition problem with a
structure similar to the classical partition problem and
3-partition. To the best of our knowledge, this is the first paper
about fractional partition. Stancu-Minasian gives a large
quantity of fractional programming problems$^{[10]}$, but all of
them are optimization problems and their constraints are not
fractional. In fact, we find MAP can be modeled as a fractional
programming problem. Not only the objective function will be
fractional, but also the constraints are fractional. Some
approximation algorithms will be presented in our subsequent papers
through modeling MAP as a fractional programming problem.

\end{document}